\documentclass{article}
\usepackage{amsmath}
\usepackage{amsfonts}
\usepackage{amssymb}
\usepackage[curve,matrix,arrow]{xy}
\def\qedsymbol{~\hfill QED.}
\def\proofname{Proof.}
\newenvironment{Proof}{\par\noindent{\it\proofname}}{{\unskip\nobreak\hfill{\it\qedsymbol}}\par\vskip 9pt}
\newenvironment{Proof*}{\par\noindent}{{\unskip\nobreak\hfill{\it\qedsymbol}}\par\vskip 9pt}
\numberwithin{equation}{section}
\newtheorem{Thm}{Theorem}[section]
\newtheorem{Lem}[Thm]{Lemma}
\newtheorem{Prop}[Thm]{Proposition}
\newtheorem{Fact}[Thm]{Fact}
\newtheorem{Cor}{Corollary}[Thm]

\newtheorem{Expl}[Thm]{Example}
\newtheorem{Rem}[Thm]{Remark}

\newcommand{\im}{\operatorname{im}}
\newcommand{\cat}{\operatorname{cat}}
\newcommand{\G}{{\Omega}}
\newcommand{\integral}{\mathbb Z}
\newcommand{\complex}{\mathbb C}
\newcommand{\quaternionic}{\mathbb H}

\newcommand{\homeo}{\approx}

\newcommand{\proj}{\operatorname{pr}}
\newcommand{\incl}{\operatorname{in}}
\newcommand{\comp}{\smash{\lower-.1ex\hbox{\scriptsize$\circ$}}}
\newcommand{\fatvee}{\operatorname{T}}
\title{
Lusternik-Schnirelmann category of a sphere-bundle over a sphere
}
\author{Norio Iwase\thanks{%
Dedicated to Professor J. R. Hubbuck on his 60th birthday.
\endgraf
This research was supported by Max-Planck-Institute f\"ur Mathematik
\endgraf
{\it MSC:\/} Primary 55M30, Secondary 55P35, 55Q25, 55R35, 55S36.
\endgraf
{\it Keywords and phrases.\/} Lusternik-Schnirelmann category, higher Hopf invariant, sphere bundle, manifold, Ganea conjecture.}\\
\\
\small {\it Address\/}:
Graduate School of Mathematics,
Kyushu University, Japan.
\\
\small {\it e-mail\/}: 
iwase@math.kyushu-u.ac.jp
}
\date{\today}
\begin{document}
\maketitle
%
%       ABSTRACT
%
\begin{abstract}\noindent%
A criterion to determine the L-S category of a total space of a sphere-bundle over a sphere is given in terms of homotopy invariants of its characteristic map, and thus providing a complete answer to Ganea's Problem 4.
As a result, we obtain a necessary and sufficient condition for such a total space $N$ to have the same L-S category as its `once punctured submanifold' $N\smallsetminus\{P\}$, $P \in N$.
Also a necessary condition for such a total space $M$ to satisfy Ganea's conjecture is obtained.
\end{abstract}
%
%       INTRODUCTION
%
\setcounter{section}{0}
\section{Introduction}\label{sect:introduction}
\par\noindent
The (normalised) L-S category $\cat(X)$ of $X$ is the least number $m$ such that there is a covering of $X$ by $m+1$ open subsets each of which is contractible in $X$, which is the least number $m$ such that the diagonal map $\Delta_{m+1} : X \to \prod^{m+1}X$ can be compressed into the `fat wedge' $\fatvee^{m+1}(X)$ (see James \cite{James:ls-category} and Whitehead \cite{Whitehead:elements}).
By definition, we have $\cat(\{\ast\}) = 0$.

This simple definition, however, does not suggest a simple way of calculation.
In fact, to {\it determine the L-S category of a sphere-bundle over a sphere in terms of homotopy invariants of its characteristic map} is listed as Problem 4 of Ganea \cite{Ganea:conjecture} in 1971.
Although a tight connection between L-S category and Bar resolution is pointed out by Ginsburg \cite{Ginsburg:ls-spectral-sequence} in 1963, this homological approach is not strong enough to solve Ganea's problems on L-S category.

Ganea's Problem 2 is also a basic problem on $\cat(X{\times}S^n)$, where we easily see $\cat(X{\times}S^n)$ = $\cat(X)$ or $\cat(X)+1$: {\it Can the latter case only occur on any $X$ and $n \geq 1$?}
The affirmative answer had become known as ``Ganea's conjecture'' or ``the Ganea conjecture'' (see James \cite{James:ls-category-survey}), particularly for manifolds.
By Singhof \cite{Singhof:minimal-cover} followed by Montejano \cite{Montejano:singhof}, G\'omez-Larra\~naga and Gonz\'alez-Acu\~na \cite{GG:ls-cat_3m} and Rudyak \cite{Rudyak:ls-cat_mfds,Rudyak:ls-cat_mfds2}, the conjecture is validated for a large class of manifolds.
%A major advance in this subject was made by Jessup \cite{Jessup:ls-cat_0} and Hess \cite{Hess:ls-cat_0} working in the rational category: the rational version of the conjecture is true for $n \geq 2$.
%Also by Singhof \cite{Singhof:minimal-cover} and later by Rudyak \cite{Rudyak:ls-cat_mfds} and \cite{Rudyak:ls-cat_mfds2}, the conjecture is true for a large class of manifolds.

The first closed manifold counter-example to the conjecture was given
by the author \cite{Iwase:counter-ls-m} as a total space of a sphere-bundle over a sphere, using concrete computations of Toda brackets depending on results by Toda \cite{Toda:composition-methods} and Oka \cite{Oka:odd-component_1}.
Also, Pascal Lambrecht, Don Stanley and Lucile Vandembroucq \cite{LSV:punctured-m} and the author \cite{Iwase:counter-ls-m} provided manifolds each of which has the same L-S category as its once punctured submanifold.
%He also knows that 
%In \cite{LSV:punctured-m}, Pascal Lambrecht, Don Stanley and Lucile
%Vandembroucq constructed another manifold with such property.
%to have the same L-S category as its once punctured submanifold.

The purpose of this paper is to determine the L-S category of a sphere-bundle over a sphere in terms of a primary homotopy invariant of the characteristic map of a bundle, providing simpler proofs of manifold examples in \cite{Iwase:counter-ls-m}.
Using it, we could obtain many closed manifolds each of which has the same L-S category as its once punctured submanifold and many closed manifold counter-examples to Ganea's conjecture on L-S category.

Throughout this paper, we follow the notations in \cite{Iwase:counter-ls,Iwase:counter-ls-m}:
In particular for a map $f : S^{k} \to X$, a homotopy set of higher Hopf invariants $H_m^S(f) = \{[H^{\sigma}_m(f)]\,\vert\,\text{$\sigma$ is a structure map of $\cat{X}{\leq}m$}\}$ (or its stabilisation ${\mathcal H}_m^S(f) = {\Sigma}^{\infty}_{\ast}H_m^S(f)$) is referred simply as {\it a (stabilised) higher Hopf invariant} of $f$, which plays a crucial role in this paper.
For $f : S^{k} \to S^{\ell}$, we identify ${H}_1^S(f)$ and ${\mathcal H}_1^S(f)$ with their unique elements, ${H}_1(f)$ and ${\mathcal H}_1(f)={\Sigma}^{\infty}{H}_1(f)$.

The author would like to express his gratitude to Hans Baues, Hans Scheerer, Daniel Tanr{\'e}, Fred Cohen, Yuli Rudyak and John Harper for valuable conversations and Max-Planck-Institut f\"ur Mathematik for its hospitality during the author's stay in Bonn.
\par
\section{L-S category of a sphere-bundle over a sphere}\label{sect:cat-sphere-bundle-over-sphere}
\par
Let $r\geq1,t\geq0$ and $E$ be a fibre bundle over $S^{t+1}$ with fibre $S^{r}$.
Then $E$ can be described as $S^{r} \cup_{\Psi} S^r{\times}D^{t+1}$, with $\Psi : S^{r} \times S^{t} \to S^{r}$ (see Whitehead \cite{Whitehead:elements}).
Hence $E$ has a CW decomposition $S^{r} \cup_{\alpha} e^{t+1} \cup_{\psi} e^{r+t+1}$ with $\alpha : S^{t} \to S^{r}$ and $\psi : S^{r+t} \to Q = S^{r} \cup_{\alpha} e^{t+1}$ given by the following formulae:
\begin{align*}&
\alpha = \Psi\vert_{\{\ast\}{\times}S^t},\quad
%\\&
%\psi=[\iota_r,\chi_{\alpha}]^r\quad 
\psi\vert_{S^{r-1}{\times}D^{t+1}} = \chi_{\alpha}{\comp}\proj_2,~
\psi\vert_{D^{r}{\times}S^{t}} = \Psi{\comp}(\omega_{r}{\times}1_{S^t}),
\end{align*}
where we denote by $\chi_{f} : (C(A),A) \to (C_f,B)$ the characteristic map for $f : A \to B$ and let $\omega_{r} = \chi_{(\ast : S^{r-1} \to \{\ast\})}$.
When $r=1$, the L-S categories of $E$ and $Q$ are studied by several authors; especially by Singhof \cite{Singhof:minimal-cover} in the case when $r=t=1$.
We summarise known results in this case.
\begin{Fact}\label{fact:cat-bundle-core0}~
Let $r=1$.
Then we have the following.
\begin{description}
\item[($t\not=0$)]~
$\cat(Q{\times}S^{n}) = 2$,
$\cat(Q) = 1$,
$\cat(E) = 2$,
$\cat(E{\times}S^{n}) = 3$.
\item[($t=1$, $\alpha=\pm1$)]~
$\cat(Q{\times}S^{n}) = 1$,
$\cat(Q) = 0$,
$\cat(E) = 1$,
$\cat(E{\times}S^{n}) = 2$.
\item[($t=1$, $\alpha=0$)]~
$\cat(Q{\times}S^{n}) = 2$,
$\cat(Q) = 1$,
$\cat(E) = 2$,
$\cat(E{\times}S^{n}) = 3$.
\item[($t=1$, $\alpha\not=0,\pm1$)]~
$\cat(Q{\times}S^{n}) = 3$,
$\cat(Q) = 2$,
$\cat(E) = 3$,
$\cat(E{\times}S^{n}) = 4$.
\item[($t>1$)]~
$\cat(Q{\times}S^{n}) = 2$,
$\cat(Q) = 1$,
$\cat(E) = 2$,
$\cat(E{\times}S^{n}) = 3$.
\end{description}
\end{Fact}
When $r>1$, we identify $H^S_1(\alpha)$ with its unique element, say $H_1(\alpha)$, since a sphere $S^k$ has the unique structure $\sigma(S^k) : S^k \to \Sigma\G{S^k}$ for $\cat(S^{k}) = 1$, $k > 1$.
We summarise the known results (due to Berstein-Hilton \cite{BH:category}) from \cite[Facts 7.1, 7.2]{Iwase:counter-ls-m}.
\begin{Fact}\label{fact:cat-bundle-core1}~
Let $r>1$.
Then we have the following.
\begin{description}
\item[($t<r$)]~
$\cat(Q{\times}S^{n}) = 2$,
$\cat(Q) = 1$,
$\cat(E) = 2$,
$\cat(E{\times}S^{n}) = 3$.
\item[($t=r$, $\alpha={\pm}1_{S^{r}}$)]~
$\cat(Q{\times}S^{n}) = 1$,
$\cat(Q) = 0$,
$\cat(E) = 1$,
$\cat(E{\times}S^{n}) = 2$.
\item[($t=r$, $\alpha\not={\pm}1_{S^{r}}$)]~
$\cat(Q{\times}S^{n}) = 2$,
$\cat(Q) = 1$,
$\cat(E) = 2$,
$\cat(E{\times}S^{n}) = 3$.
\item[($t>r$, $H_1(\alpha)=0$)]~
$\cat(Q{\times}S^{n}) = 2$,
$\cat(Q) = 1$,
$\cat(E) = 2$,
$\cat(E{\times}S^{n}) = 3$.
\item[($t>r$, $H_1(\alpha)\not=0$)]~
$\cat(Q{\times}S^{n}) =$ $3$ or $2$,
$\cat(Q) = 2$,
$\cat(E) =$ $2$ or $3$,
$\cat(E{\times}S^{n}) =$ $3$ or $4$.
%$2\leq\cat(Q{\times}S^{n})\leq3$,
%$\cat(Q) = 2$,
%$2\leq\cat(E)\leq3$,
%$3\leq\cat(E{\times}S^{n})\leq4$.
\end{description}
\end{Fact}
By \cite{Iwase:counter-ls} and \cite[Theorem 5.2, 5.3, 7.3]{Iwase:counter-ls-m}, the following is also known.
\begin{Fact}\label{fact:cat-bundle1}~
When $r > 1$, $t \geq r$ and $\alpha\not=\pm1$, 
%In the last case in Fact \ref{fact:cat-bundle-core1}, 
we also have the following.
\begin{enumerate}
\item[(1)]~
${\Sigma}^{n}H_1(\alpha)=0$ implies $\cat(Q{\times}S^n)$ $=$ $2$, and 
${\Sigma}^{n+1}H_1(\alpha)\not=0$ implies $\cat(Q{\times}S^n)$ $=$ $3$.
\item[(2)]~
$\cat(E)$ $=$ $2$ if and only if $H^S_2(\psi)$ $\ni$ $0$, and 
$\cat(E)$ $=$ $2$ implies $\cat(E{\times}S^n)$ $=$ $3$ for all $n$.
\item[(3)]~
${\Sigma}^n_{\ast}H^S_2(\psi)$ $\ni$ $0$ implies $\cat(E{\times}S^n)$ $=$ $3$, and 
${\Sigma}^{n+r+1}h_2(\alpha)$ $\not=$ $0$ implies $\cat(E{\times}S^n)$ $=$ $4$.
\end{enumerate}
\end{Fact}
\begin{Rem}\label{rem:cat-bundle_sphere}~
When $\alpha$ is in meta-stable range, $H_1(\alpha) : S^{t} \to \G{S^{r}}{\ast}\G{S^{r}}$ is given by the second James-Hopf invariant $h_2(\alpha) : S^{t} \to \Sigma{S^{r-1}{\wedge}S^{r-1}}$ composed with an appropriate inclusion to a wedge-summand.
Thus we may regard $h_2(\alpha)=H_1(\alpha)$ when $\alpha$ is in meta-stable range.
\end{Rem}

But a higher Hopf invariant $H^S_2(\psi)$ is not very easy to determine.
% (see Section 8 of \cite{Iwase:counter-ls-m}).
Our result is as follows:
\begin{Thm}\label{lem:crucial}~
Let $\cat(Q)=2$ with $t > r > 1$, Then 
$H^S_2(\psi)$ contains $0$ if and only if ${\Sigma}^{r}H_1(\alpha) = 0$.
%\end{Thm}
More generally %we have the following result
%\begin{Thm}
%Let $E = S^r \cup_{\alpha} e^{t+1} \cup_{\psi} e^{r+t+1}$ be an $S^{r}$-bundle over $S^{t+1}$, $Q = S^r \cup_{\alpha} e^{t+1}$ and 
for a co-H-map $\beta : S^{v} \to S^{r+t}$, $H^{S}_{2}(\psi{\comp}\beta) = \beta^{\ast}H^{S}_{2}(\psi)$ contains $0$ if and only if ${\Sigma}^rH_1(\alpha){\comp}\beta = 0$.
\end{Thm}
The result is obtained by the following lemma for $Q$ of $\cat(Q)=2$ with $t > r > 1$.
\begin{Lem}\label{prop:standard_structure}~
%Let $\cat(Q)=2$ with $t > r > 1$.
%The higher Hopf invariant 
%$H^S_2(\psi)$ contains ${\Sigma}^{r}H_1(\alpha)$ or $-{\Sigma}^{r}H_1(\alpha)$.
$H^S_2(\psi) \ni \pm [(\hat{i}{\ast}1_{\G{Q}{\ast}\G{Q}}){\comp}{\Sigma}^{r}H_1(\alpha)]$, where the bottom-cell inclusion $\hat{i} : S^{r-1} \hookrightarrow \G{Q}$ denotes the adjoint of the inclusion $i : S^{r} \hookrightarrow Q$.
\end{Lem}
%\begin{Rem}
%For dimensional reasons, the image of $(\hat{i}{\ast}1_{\G{Q}{\ast}\G{Q}}){\comp}{\Sigma}^{r}H_1(\alpha)$ lies in $\G{S^r}{\ast}\G{S^r}{\ast}\G{S^r}$ in which the subspace ${S^{r-1}}{\ast}\G{S^r}{\ast}\G{S^r}$ is a direct wedge-sum summand.
%Thus $(\hat{i}{\ast}1_{\G{Q}{\ast}\G{Q}}){\comp}{\Sigma}^{r}H_1(\alpha)$ is non-trivial if and only if so is ${\Sigma}^{r}H_1(\alpha)$.
%\end{Rem}

By combining above facts with Theorem \ref{lem:crucial}, we obtain an answer to Ganea's Problem 4:
\begin{Thm}\label{thm:main-theorem}~\par
\begin{center}
\begin{tabular}{|c|c|c||c|c|c|c|}
\hline
\multicolumn{3}{|c||}{Conditions} & \multicolumn{4}{c|}{L-S category}
\\\hline
$r$& $t$& $\alpha$& ${Q{\times}S^{n}}$& ${Q}$& ${E}$& ${E{\times}S^{n}}$
\\\hline
\raisebox{-5.7ex}[0cm][0cm]{$r=1$}	%
	& $t=0$	&			& $2$	& $1$	& $2$	& $3$
\\\cline{2-7}
	&\raisebox{-2.9ex}[0cm][0cm]{$t=1$}	%
		& $\alpha=\pm1$		& $1$	& $0$	& $1$	& $2$
\\\cline{3-7}
	& 	& $\alpha=0$		& $2$	& $1$	& $2$	& $3$
\\\cline{3-7}
	& 	& $\alpha\not=0,\pm1$	& $3$	& $2$	& $3$	& $4$
\\\cline{2-7}
	& $t>1$	& 			& $2$	& $1$	& $2$	& $3$
\\\hline
\raisebox{-8.1ex}[0cm][0cm]{$r>1$}	%
	& $t<r$	& 			& $2$	& $1$	& $2$	& $3$
\\\cline{2-7}
	&\raisebox{-1.5ex}[0cm][0cm]{$t=r$}	%
		& $\alpha=\pm1$		& $1$	& $0$	& $1$	& $2$
\\\cline{3-7}
	& 	& $\alpha\not=\pm1$		& $2$	& $1$	& $2$	& $3$
\\\cline{2-7}
	&\raisebox{-4.0ex}[0cm][0cm]{$t>r$}	%
		& $H_1(\alpha)=0$	& $2$	& $1$	& $2$	& $3$
\\\cline{3-7}
	& 	& $H_1(\alpha)\not=0 ~\,\&~{\Sigma}^rH_1(\alpha)=0$%
					&\raisebox{-3.44ex}[0cm][0cm]{\begin{tabular}{c}$3$ or $2$\\[-1mm]{\small(1)}\end{tabular}}%
						&\raisebox{-2.5ex}[0cm][0cm]{$2$}%
							& $2$	& $3$
\\\cline{3-3}\cline{6-7}
	& 	& ${\Sigma}^rH_1(\alpha)\not=0$	
					&	& 	& $3$	&\begin{tabular}{c}$3$ or $4$\\[-1mm]{\small(2)}\end{tabular}
\\\hline
\end{tabular}
\end{center}
\par
{\scriptsize
(1): 
$\begin{displaystyle}
\left\{\begin{array}{l} 
\text{${\Sigma}^{n}H_1(\alpha)=0$ implies $\cat(Q{\times}S^n)=2$ and }
\\
\text{${\Sigma}^{n+1}H_1(\alpha)\not=0$ implies $\cat(Q{\times}S^n)=3$.}
\end{array}\right.
\end{displaystyle}$
~(2): 
$\begin{displaystyle}
\left\{\begin{array}{l}
\text{${\Sigma}^{r+n}H_1(\alpha)=0$ implies $\cat(E{\times}S^n)=3$ and }
\\
\text{${\Sigma}^{r+n+1}h_2(\alpha)\not=0$ implies $\cat(E{\times}S^n)=4$.}
\end{array}\right.
\end{displaystyle}$
}
\end{Thm}
%
%	EXAMPLES
%
\section{Applications and examples}\label{sect:examples}
\par\noindent
Firstly, Theorem \ref{thm:main-theorem} yields the following result.
\begin{Thm}\label{thm:punctured-m}~
Let a manifold $N$ be the total space of a $S^{r}$-bundle over $S^{t+1}$ with a characteristic map $\Psi : S^{r}{\times}S^{t} \to S^{r}$, $t > r > 1$, and let $\alpha=\Psi\vert_{S^{t}}$.
Then $\cat(N\smallsetminus\{P\})=\cat(N)$ if and only if $H_1(\alpha)\not=0$ and ${\Sigma}^{r}H_1(\alpha)=0$.
%The converse is also true.
\end{Thm}
This theorem provides the following examples.
\begin{Expl}~
Let $p$ be an odd prime and $\alpha=\eta_2{\comp}\alpha_1(3){\comp}\alpha_1(2p)$.
Then we have that $H_1(\alpha)=\alpha_1(3){\comp}\alpha_1(2p)\ne0$ and ${\Sigma}^2H_1(\alpha)=0$ by \cite{Toda:composition-methods}.
Let $N_p \to S^{4p-2}$ be the bundle with fibre $S^{2}$ induced by $\Sigma(\alpha_1(3){\comp}\alpha_1(2p)) : S^{4p-2} \to S^4$ from the bundle ${\complex}P^3 \to {\quaternionic}P^1=S^4$ with fibre $Sp(1)/U(1)=S^2$.
By the argument given in \cite{Iwase:counter-ls-m} shows that $N_p$ has a CW-decomposition as $N_p \homeo S^2 \cup_{\alpha} e^{4p-2} \cup_{\psi} e^{4p}$.
Then Theorem \ref{thm:punctured-m} implies that $\cat(N_p)=\cat(N_p\smallsetminus\{P\})=2$.
\end{Expl}
\begin{Expl}[\cite{Iwase:counter-ls-m}]~
Let $p$ be a prime $\geq 5$ and $\alpha=\eta_2{\comp}\alpha_1(3){\comp}\alpha_2(2p)$ as in \cite{Iwase:counter-ls-m}.
Then we have that $H_1(\alpha)=\alpha_1(3){\comp}\alpha_2(2p)\ne0$ and ${\Sigma}^2H_1(\alpha)=0$ by \cite{Toda:composition-methods}.
Let $L_p \to S^{6p-4}$ be the bundle with fibre $S^{2}$ induced by $\Sigma(\alpha_1(3){\comp}\alpha_2(2p)) : S^{6p-4} \to S^4$ from the bundle ${\complex}P^3 \to {\quaternionic}P^1=S^4$ with fibre $Sp(1)/U(1)=S^2$.
By the argument given in \cite{Iwase:counter-ls-m} shows that $L_p$ has a CW-decomposition as $L_p \homeo S^2 \cup_{\alpha} e^{6p-4} \cup_{\psi} e^{6p-2}$.
Then Theorem \ref{thm:punctured-m} implies that $\cat(L_p)=\cat(L_p\smallsetminus\{P\})=2$.
\end{Expl}
\par
Secondly, Theorem \ref{thm:main-theorem} also yields the following result.
\begin{Thm}\label{thm:counter-ganea-m}~
Let a manifold $M$ be the total space of a $S^{r}$-bundle over $S^{t+1}$ with a characteristic map $\Psi : S^{r}{\times}S^{t} \to S^{r}$, $t > r > 1$, and let $\alpha=\Psi\vert_{S^{t}}$.
If ${\Sigma}^{r}H_1(\alpha) \ne 0$ and ${\mathcal H}_1(\alpha) = 0$, then $M$ is a counter-example to the Ganea's conjecture on L-S category; more precisely, $\cat(M)=\cat(M{\times}S^n)=3$ if ${\Sigma}^{r}H_1(\alpha) \ne 0$ and ${\Sigma}^{n+r}H_1(\alpha) = 0$.
\end{Thm}
This theorem provides the following  manifold counter examples to Ganea's conjecture on L-S category.
\begin{Expl}~
Let $p=2$ and $\alpha=\eta_2{\comp}\eta^2_3{\comp}\epsilon_5$.
Then we have that $H_1(\alpha)=\eta^2_3{\comp}\epsilon_5\ne0$, ${\Sigma}^2H_1(\alpha)\ne0$ and ${\Sigma}^6H_1(\alpha)=0$ by \cite{Toda:composition-methods}.
Let $M_2 \to S^{14}$ be the bundle with fibre $S^{2}$ induced by $\Sigma(\eta^2_3{\comp}\epsilon_5) : S^{14} \to S^4$ from the bundle ${\complex}P^3 \to {\quaternionic}P^1=S^4$ with fibre $Sp(1)/U(1)=S^2$.
By the argument given in \cite{Iwase:counter-ls-m} shows that $M_2$ has a CW-decomposition as $M_2 \homeo S^2 \cup_{\alpha} e^{14} \cup_{\psi} e^{16}$.
Then Theorem \ref{thm:counter-ganea-m} implies that $\cat(M_2{\times}S^{n})=\cat(M_2)=3$ for $n \geq 4$.
\end{Expl}
\begin{Expl}[\cite{Iwase:counter-ls-m}]~
Let $p=3$ and $\alpha=\eta_2{\comp}\alpha_1(3){\comp}\alpha_2(6)$ as in \cite{Iwase:counter-ls-m}.
Then we have that $H_1(\alpha)=\alpha_1(3){\comp}\alpha_2(6)\ne0$, ${\Sigma}^2H_1(\alpha)\ne0$ and ${\Sigma}^4H_1(\alpha)=0$ by \cite{Toda:composition-methods}.
Let $M_3 \to S^{14}$ be the bundle with fibre $S^{2}$ induced by $\Sigma(\alpha_1(3){\comp}\alpha_2(6)) : S^{14} \to S^4$ from the bundle ${\complex}P^3 \to {\quaternionic}P^1=S^4$ with fibre $S^2$.
By the argument given in \cite{Iwase:counter-ls-m} shows that $M_3$ has a CW-decomposition as $M_3 \homeo S^2 \cup_{\alpha} e^{14} \cup_{\psi} e^{16}$.
Then Theorem \ref{thm:counter-ganea-m} implies that $\cat(M_3{\times}S^{n})=\cat(M_3)=3$ for $n \geq 2$.
\end{Expl}
Finally, Theorem \ref{lem:crucial} and \cite[Theorem 5.2]{Iwase:counter-ls-m} imply the following result.
\begin{Thm}
Let a manifold $X$ be the total space of a $S^{r}$-bundle over $S^{t+1}$ with a characteristic map $\Psi : S^{r}{\times}S^{t} \to S^{r}$, $t > r > 1$, and let $\alpha=\Psi\vert_{S^{t}}$.
When $H_1(\alpha) \ne 0$ and $\beta$ is a co-H-map, we obtain that $X(\beta) = S^{r} \cup_{\alpha} e^{t+1} \cup_{\psi{\comp}\beta} e^{v+1}$ is of $\cat(X(\beta))=3$ if and only if ${\Sigma}^{r}H_1(\alpha){\comp}\beta \not= 0$. 
\end{Thm}
%
%	CRUCIAL LEMMA
%
\section{Proof of Lemma \ref{prop:standard_structure}}\label{section:standard_structure}
\par\noindent
Let $\cat(Q)=2$ with $t > r > 1$.
In the remainder of this paper, we distinguish a map from its homotopy class to make the arguments clear.
\par
Here, let us recall the definition of a {\it relative Whitehead product}:
For maps $f : {\Sigma}X \to M$ and $g : (C(Y),Y) \to (K,L)$, we denote by $[f,g]^{\rm rel} : X{\ast}Y=C(X){\times}Y \cup X{\times}C(Y) \to M{\times}L \cup \{\ast\}{\times}K$ the relative Whitehead product, which is given by
\begin{align*}&
[f,g]^{\text{\rm rel}}\vert_{C(X){\times}Y}(t{\wedge}x,y) = (f(t{\wedge}x),g(y))\quad\text{and}
%\quad
\\&
[f,g]^{\text{\rm rel}}\vert_{X{\times}C(Y)}(x,t{\wedge}y) = (\ast,g(t{\wedge}y)).
\end{align*}
Also a pairing $F : M{\times}L \to M$ with axes $1_M$ and $h : L \to M$ (see Oda \cite{Oda:pairing}) determines a map 
\begin{equation*}
(F \cup \chi_h) : (M{\times}L \cup \{\ast\}{\times}K) \to (M \cup_hK,M)
\end{equation*}
by $(F \cup \chi_h)\vert_{ M{\times}L} = F$ and $(F \cup \chi_h)\vert_{\{\ast\}{\times}K} = \chi_h$, where $\chi_h : (K,L) \to (M\cup_hK,M)$ is a relative homeomorphism given by the restriction of the identification map $M \cup K \to M\cup_hK$.
Then we can easily see that $\psi : S^{r+t} \to Q$ is given as 
\begin{equation}\label{eq:psi}
\psi = (\Psi\cup\chi_{\alpha}){\comp}[\iota_r,C(\iota_{t})],
\end{equation}
where $\iota_k : S^k \to S^k$  and $C(\iota_{k}) : C(S^k) \to C(S^k)$ denote the identity maps.
\par
We denote by $j^{Q}_{i} : P^{i}(\G{Q}) \overset{j^Q_i}\hookrightarrow P^{\infty}(\G{Q})$ the classifying map of the fibration $p^{\G{Q}}_{i} : E^{i+1}(\G{Q}) \to P^{i}(\G{Q})$ and $e^{Q}_{i} = e^{Q}_{\infty}{\comp}j^{Q}_{i}$, where $e^{Q}_{\infty} : P^{\infty}(\G{Q}) \to Q$ is a homotopy equivalence extending the evaluation map $e^Q_1=ev : {\Sigma}\G{Q} \to Q$.
Let $\sigma_{\infty}$ be the homotopy inverse of $e^Q_{\infty}$.
Then we may assume that $\sigma_{\infty}\vert_{S^r}=j^Q_1{\comp}\sigma(S^r)$ for dimensional reasons.
\begin{Prop}\label{prop:pull-back}~
The following is a commutative diagram where the lower squares are pull-back diagrams.
\begin{equation}\label{eq:pushpull_3}
\begin{xy}
\xymatrix{
E^3(\G{Q})
	\ar[dd]_{p^{\G{Q}}_2}
	\ar@{=}[rr]
&
&
\G{Q}{\ast}E^2(\G{Q})
	\ar[dd]_{[j^Q_1,j^Q_2{\comp}\chi_{p^{\G{Q}}_1}]^{\text{\rm rel}}}
	\ar@{=}[rr]
&
&
\G{Q}{\ast}\G{Q}{\ast}\G{Q}
	\ar[dd]_{[e^Q_1,(e^Q_1{\times}e^Q_1){\comp}{\chi}_{[\iota,\iota]}]^{\text{\rm rel}}}
\\
\\
P^2(\G{Q})
	\ar[dd]^{e^Q_2}
	\ar[rr]_{\hat{\Delta}_{Q}\qquad\qquad}
&
&
P^{\infty}(\G{Q}){\times}{\Sigma\G{Q}} \cup \{\ast\}{\times}P^{\infty}(\G{Q})
	\ar[dd]^{e^Q_{\infty}{\times}e^{Q}_1 \cup {\ast}{\times}e^Q_{\infty}}
	\ar[rr]
&
&
\fatvee^3{Q}
	\ar@{^{(}->}[dd]
\\
\\
Q
	\ar@{.>}[uurr]^{\sigma'_0}
	\ar@{.>}@/^/[uu]^{\sigma_0}
	\ar[rr]_{\Delta_{Q}}
&
&
Q{\times}Q
	\ar[rr]_{1_Q{\times}\Delta_{Q}}
&
&
Q{\times}Q{\times}Q.
}
\end{xy}
\end{equation}
\end{Prop}
\begin{Rem}~
The homotopy fibre $\G{Q}{\ast}\G{Q}{\ast}\G{Q} \to \fatvee^3{Q}$ of the inclusion 
$$\fatvee^3{Q} = Q{\times}(Q{\vee}Q) \cup \{\ast\}{\times}(Q{\times}Q) \hookrightarrow Q{\times}(Q{\times}Q)$$ 
is given by a relative Whitehead product $[e^Q_1,(e^Q_1{\times}e^Q_1){\comp}{\chi}_{[\iota,\iota]}]^{\text{\rm rel}}% = \fatvee^3(e^Q_1){\comp}[\iota,{\chi}_{[\iota,\iota]}]^{\text{\rm rel}}
$, where $\iota$ denotes the identity $1_{{\Sigma}\G{Q}}$ and 
$${\chi}_{[\iota,\iota]} : (C(\G{Q}{\ast}\G{Q}),\G{Q}{\ast}\G{Q}) \to ({\Sigma}\G{Q}{\times}{\Sigma}\G{Q},{\Sigma}\G{Q}{\vee}{\Sigma}\G{Q})$$ 
denotes a relative homeomorphism.
\end{Rem}
A lifting $\sigma'_0$ of $\Delta_{Q}$ in diagram (\ref{eq:pushpull_3}) is given by the following data:
\begin{align*}&
\sigma'_0\vert_{S^r} = ((j^Q_1{\comp}\sigma(S^r))\times\sigma(S^r)){\comp}\Delta_{S^r}%\quad \text{and, for $u{\wedge}x \in (0,1]{\times}S^t/\{1\}{\times}S^r \subset Q$,}
\intertext{and $\sigma'_0\vert_{Q{\smallsetminus}S^r}$ for $u{\wedge}x \in (0,1]{\times}S^t/\{1\}{\times}S^r \subset Q$,}&
\sigma'_0(u{\wedge}x) = \left\{
\begin{array}{ll}
((j^Q_1{\comp}\sigma(S^r)){\comp}\alpha\times\sigma(S^r){\comp}\alpha){\comp}H_t(2u{\wedge}x)&\text{if $u \leq \frac{1}{2}$}
\\
(\hat{\chi}_{\alpha}(2u-1,x_1),\hat{\chi}_{\alpha}(2u-1,x_2)),~ \mu_t(x)=(x_1,x_2)&\text{if $u \geq \frac{1}{2}$},
\end{array}\right.
\end{align*}
where $H_t$ is a homotopy $\Delta_{S^t} \sim \mu_{t}$ in $S^t{\times}S^t$, $\mu_k={\Sigma}^{k-1}\mu_{1} : S^k \to S^k \vee S^k$ denotes the unique co-H-structure of $S^k$ and $\hat{\chi}_{\alpha}$ is a null-homotopy $\sigma_{\infty}{\comp}\chi_{\alpha} : (C(S^{t}),S^t) \to (Q,S^r) \to (P^{\infty}(\G{Q}),\im(j^Q_1{\comp}\sigma(S^r)))$ of $j^Q_1{\comp}\sigma(S^r){\comp}\alpha \sim \ast$.
\par
Since the lower left square of diagram (\ref{eq:pushpull_3}) is a homotopy pullback diagram, $\sigma'_0$ and the identity $1_Q$ defines a lifting $\sigma_0 : Q \to P^2(\G{Q})$ of $1_Q$.
\begin{Proof*}{\it Proof of Proposition \ref{eq:pushpull_3}.}
By \cite[Lemma 2.1]{Iwase:counter-ls} with $(X,A) = (P^{\infty}(\G{Q}),\{\ast\})$, $(Y,B) = (P^{\infty}(\G{Q}),\Sigma\G{Q})$, $Z = P^{\infty}(\G{Q})$ and $f = g = 1_{P^{\infty}(\G{Q})}$, we have the following homotopy pushout-pullback diagram:
\begin{equation}\label{eq:pushpull_1}
\begin{xy}
\xymatrix{
E^{2}(\G{Q})
	\ar[rrr]
	\ar[dd]_{p^{\G{Q}}_1}
&
&
&
\{\ast\}
	\ar[dd]
&
\\
	\ar@{{}}^{\qquad\qquad HPO}
&
&
&
&
\\
{\Sigma}\G{Q}
	\ar[rrr]
&
&
&
P^{2}(\G{Q})
	\ar[r]^{\hat{\Delta}_{Q}\qquad\qquad\quad}
	\ar[dd]_{e^{Q}_2}
&
P^{\infty}(\G{Q}){\times}{\Sigma\G{Q}} \cup \{\ast\}{\times}P^{\infty}(\G{Q})
	\ar[dd]^{e^{Q}_{\infty}{\times}e^{Q}_1 \cup {\ast}{\times}e^{Q}_{\infty}}
\\
&
&
&
	\ar@{{}}^{\qquad\qquad HPB}
&
\\
&
&
&
Q
	\ar[r]_{\Delta_{Q}\qquad}
&
Q{\times}Q,
}
\end{xy}
\end{equation}
where we replaced $P^{\infty}(\G{Q})$ by $Q$ in the bottom, since $P^{\infty}(\G{Q})$ is the homotopy equivalent with $Q$ by $e^{Q}_{\infty} : P^{\infty}(\G{Q}) \to Q$ and $\sigma_{\infty} : Q \to P^{\infty}(\G{Q})$.
\par
By \cite[Lemma 2.1]{Iwase:counter-ls} with $(X,A) = (P^{\infty}(\G{Q}),\{\ast\})$, $(Y,B) = (P^{\infty}(\G{Q}),\Sigma\G{Q})$, $Z=\{\ast\}$ and $f=g=\ast$, we have the following homotopy pushout-pullback diagram:
\begin{equation}\label{eq:pushpull_2}
\begin{xy}
\xymatrix{
\G{Q}{\times}E^2(\G{Q})
	\ar[rr]_{\proj_1}
	\ar[dd]_{\proj_2}
&
&
\G{Q}
	\ar[dd]
&
\\
	\ar@{{}}^{\qquad\qquad HPO}
&
&
\\
E^2(\G{Q})
	\ar[rr]
&
&
\G{Q}{\ast}E^2(\G{Q})
	\ar[r]^{[j^Q_1,j^Q_2{\comp}\chi_{p^{\G{Q}}_1}]^{\text{\rm rel}}\qquad\qquad}
	\ar[dd]_{}
&
P^{\infty}(\G{Q}){\times}{\Sigma\G{Q}}\cup\{\ast\}{\times}P^{\infty}(\G{Q})
	\ar[dd]^{e^{Q}_{\infty}{\times}e^{Q}_1\cup{\ast}{\times}e^{Q}_{\infty}}
\\
&
&
	\ar@{{}}^{\qquad\qquad HPB}
&
\\
&
&
\{\ast\}
	\ar[r]_{}
&
Q{\times}Q,
}
\end{xy}
\end{equation}
where $\chi_{p^{\G{Q}}_1} : (C(E^2(\G{Q})),E^2(\G{Q})) \to (P^{2}(\G{Q}),\Sigma\G{Q})$ is a relative homeomorphism.
\par
The above constructions give a standard $\G{Q}$-projective plane $P^{2}(\G{Q})$ and a standard projection $p^{\G{Q}}_2 : E^3(\G{Q}) \to P^{2}(\G{Q})$.
In fact, the diagonal map $\Delta^3_{Q} : Q \to Q{\times}Q{\times}Q$ is the composition $(1_Q{\times}\Delta_{Q}){\comp}\Delta_{Q}$ and there is the following homotopy pushout-pullback diagram by \cite[Lemma 2.1]{Iwase:counter-ls} with $(X,A) = (Q,\{\ast\})$, $(Y,B) = (Q{\times}Q,Q{\vee}Q)$, $Z=Q{\times}Q$, $f=\proj_1$ and $g=\Delta_{Q}{\comp}\proj_2$:
\begin{equation*}\label{eq:pushpull_0}
\begin{xy}
\xymatrix{
\{\ast\}{\times}{\Sigma}\G{Q}
	\ar@{^{(}->}[r]
	\ar@{^{(}->}[dd]
&
\{\ast\}{\times}{Q}
	\ar[dd]
&
\\
	\ar@{{}}^{\qquad\qquad HPO}
&
&
\\
Q{\times}{\Sigma}\G{Q}
	\ar[r]
&
P^{\infty}(\G{Q}){\times}{\Sigma\G{Q}} \cup \{\ast\}{\times}P^{\infty}(\G{Q})
	\ar[r]
	\ar[dd]_{e^Q_{\infty}{\times}e^{Q}_1 \cup {\ast}{\times}e^Q_{\infty}}
&
\fatvee^3{Q}
	\ar@{^{(}->}[dd]
\\
&
	\ar@{{}}^{\qquad\qquad HPB}
&
\\
&
Q{\times}Q
	\ar[r]_{1_{Q}{\times}\Delta_{Q}}
&
Q{\times}Q{\times}Q.
}
\end{xy}
\end{equation*}
%Thus (\ref{eq:pushpull_1}) and (\ref{eq:pushpull_2}) yield the following commutative diagram:
By combining this diagram with diagrams (\ref{eq:pushpull_1}) and (\ref{eq:pushpull_2}), we obtain the desired diagram.
\end{Proof*}
%a formula $\sigma'_0\vert_{S^r} = (\sigma(S^r)\times\sigma(S^r)){\comp}\Delta_{S^r}$ together with two homotopies $(\sigma(S^r){\comp}\alpha\times\sigma(S^r){\comp}\alpha){\comp}H_t$ and $\hat{\chi}_{\alpha}{\vee}\hat{\chi}_{\alpha}$, 
%for any element $u{\wedge}x \in (0,1]{\times}S^t/\{1\}{\times}S^r \subset Q$, 
%Then by using the homotopies $(\sigma(S^r){\comp}\alpha\times\sigma(S^r){\comp}\alpha){\comp}H_t$ and $H_r : \Delta_{S^r} \sim \mu_r$, we obtain the following.
Since there is a right action of $S^t{\times}S^t$ on $S^r{\times}S^r$ by $\Psi^2=(\Psi{\times}\Psi){\comp}(1{\times}T{\times}1) : S^r{\times}S^r{\times}S^t{\times}S^t \to S^r{\times}S^r$, we obtain the following.
\begin{Prop}\label{Prop:sigma}~
The map $\sigma'_0{\comp}\psi : S^{t} \to P^{\infty}(\G{Q}){\times}{\Sigma}\G{Q} \cup \{\ast\}{\times}P^{\infty}(\G{Q})$ satisfies 
\begin{align*}
\sigma'_0{\comp}\psi 
&\sim (((j^Q_1{\comp}\sigma(S^r))\times\sigma(S^r)){\comp}\Psi^2_0 \cup (\hat{\chi}_{\alpha}\vee\hat{\chi}_{\alpha})){\comp}[\mu_{r},C(\mu_{t})]^{\text{\rm rel}},
\end{align*}
where $\Psi^2_{0} = \Psi^2\vert_{(S^{r}{\vee}S^{r}){\times}(S^{t}{\vee}S^{t})} : (S^{r}{\vee}S^{r}){\times}(S^{t}{\vee}S^{t}) \to S^{r}{\vee}S^{r}$.
\end{Prop}
\begin{Proof}~
By (\ref{eq:psi}), we know $\sigma'_0{\comp}\psi = \sigma'_0{\comp}(\Psi\cup\chi_{\alpha}){\comp}[\iota_r,C(\iota_t)]^{\text{\rm rel}} = \sigma'_0{\comp}(\Psi\cup\chi_{\alpha}) = (\sigma'_0\vert_{\im\sigma(S^r)}{\comp}\Psi\cup\sigma'_0{\comp}\chi_{\alpha}){\comp}[\iota_r,C(\iota_t)]^{\text{\rm rel}}$, where we have
\begin{align*}&
\sigma'_0\vert_{\im\sigma(S^r)}{\comp}\Psi = j^Q_1{\comp}\sigma(S^r){\comp}\Delta_{S^r}{\comp}\Psi = j^Q_1{\comp}\sigma(S^r){\comp}\Psi^2{\comp}(\Delta_{S^r}{\times}\Delta_{S^t})\quad\text{and}
\\&
\sigma'_0{\comp}\chi_{\alpha} = ((j^Q_1{\comp}\sigma(S^r)){\comp}\alpha{\times}(j^Q_1{\comp}\sigma(S^r)){\comp}\alpha){\comp}H_{t} + (\hat{\chi}_{\alpha}{\vee}\hat{\chi}_{\alpha}){\comp}C(\mu_t),
\end{align*}
where the addition denotes the composition of homotopies.
Using the same homotopy $H_t : \Delta_{S^t} \sim \mu_t$, we obtain homotopies
\begin{align*}&
\sigma'_0\vert_{\im\sigma(S^r)}{\comp}\Psi \sim j^Q_1{\comp}\sigma(S^{r}){\comp}\Psi^2{\comp}(\Delta_{S^{r}}{\times}\mu_t),
\quad\text{and}\quad
%\\&
\sigma'_0{\comp}\chi_{\alpha} \sim (\hat{\chi}_{\alpha}{\vee}\hat{\chi}_{\alpha}){\comp}C(\mu_t)
\end{align*}
which fit together into a homotopy
\begin{equation*}
\sigma'_0{\comp}(\Psi\cup\chi_{\alpha}) \sim (((j^Q_1{\comp}\sigma(S^r))\times\sigma(S^r)){\comp}\Psi^2{\comp}(\Delta_{S^r}{\times}\mu_t) \cup (\hat{\chi}_{\alpha}\vee\hat{\chi}_{\alpha}){\comp}C(\mu_t)).
\end{equation*}
Then the homotopy $H_r : \Delta_{S^r} \sim \mu_{r}$ gives the homotopy relation
\begin{align*}
\sigma'_0{\comp}\psi &\sim (((j^Q_1{\comp}\sigma(S^r))\times\sigma(S^r)){\comp}\Psi^2{\comp}(\mu_{r}{\times}\mu_t) \cup (\hat{\chi}_{\alpha}\vee\hat{\chi}_{\alpha}){\comp}C(\mu_t)){\comp}[\iota_r,C(\iota_t)]^{\text{\rm rel}},
\end{align*}
%\\&\quad
which yields $\sigma'_0{\comp}\psi \sim (((j^Q_1{\comp}\sigma(S^r))\times\sigma(S^r)){\comp}\Psi^2_0 \cup (\hat{\chi}_{\alpha}\vee\hat{\chi}_{\alpha})){\comp}[\mu_{r},C(\mu_{t})]^{\text{\rm rel}}$.
\end{Proof}
Hence by the definition of $\sigma_0$ and $\psi$, we obtain the following.
\begin{Prop}\label{prop:relative-whitehead}~
$\hat{\Delta}_{Q}{\comp}p^{\G{Q}}_2{\comp}H^{\sigma_{0}}_2(\psi) \sim [j^Q_1{\comp}\sigma(S^{r}),\hat{\chi}_{\alpha}]^{\text{\rm rel}}$, 
\end{Prop}
\begin{Proof}~
By the definition of $\sigma_0$, we obtain 
\begin{equation*}
\hat{\Delta}_{Q}{\comp}\sigma_0{\comp}\psi \sim \sigma'_0{\comp}\psi \sim (((j^Q_1{\comp}\sigma(S^r))\times\sigma(S^r)){\comp}\Psi^2_0 \cup (\hat{\chi}_{\alpha}\vee\hat{\chi}_{\alpha})){\comp}[\mu_{r},C(\mu_{t})]^{\text{\rm rel}}.
\end{equation*}
Let $\incl_i : Z \to Z{\vee}Z$ be the inclusion to the $i$-th factor.
Then $[\mu_{r},C(\mu_{t})]^{\text{\rm rel}} : S^{r+t} \to (S^r{\vee}S^r){\times}(S^t{\vee}S^t)$ can be deformed as 
\begin{align*}&
[\mu_{r},C(\mu_{t})]^{\text{\rm rel}} \sim [\incl_1{\comp}\iota_{r}+\incl_2{\comp}\iota_{r},\incl_1{\comp}C(\iota_{t})+\incl_2{\comp}C(\iota_{t})]^{\text{\rm rel}}
\\&\quad
\sim [\incl_1{\comp}\iota_r,\incl_1{\comp}C(\iota_{t})]^{\text{\rm rel}} + [\incl_2{\comp}\iota_r,\incl_2{\comp}C(\iota_{t})]^{\text{\rm rel}}
\\&\quad\qquad\qquad
+ [\incl_2{\comp}\iota_r,\incl_1{\comp}C(\iota_{t})]^{\text{\rm rel}} + [\incl_1{\comp}\iota_r,\incl_2{\comp}C(\iota_{t})]^{\text{\rm rel}}
\\&\quad
\sim [\incl_1{\comp}\iota_r,\incl_1{\comp}C(\iota_{t})]^{\text{\rm rel}} + [\incl_2{\comp}\iota_r,\incl_2{\comp}C(\iota_{t})]^{\text{\rm rel}}
\\&\quad\qquad\qquad
+ [\incl_2{\comp}\iota_r,\incl_1{\comp}C(\iota_{t})]^{\text{\rm rel}} + [\incl_1{\comp}\iota_r,\incl_2{\comp}C(\iota_{t})]^{\text{\rm rel}}
\end{align*}
in $(S^r{\vee}S^r){\times}(S^t{\vee}S^t)$.
Thus we have 
\begin{align*}&
\hat{\Delta}_{Q}{\comp}\sigma_0{\comp}\psi \sim (((j^Q_1{\comp}\sigma(S^r))\times\sigma(S^r)){\comp}\Psi^2_0 \cup (\hat{\chi}_{\alpha}\vee\hat{\chi}_{\alpha})){\comp}[\incl_1{\comp}\iota_r,\incl_1{\comp}C(\iota_{t})]^{\text{\rm rel}}
\\&\qquad\qquad\qquad
+ (((j^Q_1{\comp}\sigma(S^r))\times\sigma(S^r)){\comp}\Psi^2_0 \cup (\hat{\chi}_{\alpha}\vee\hat{\chi}_{\alpha})){\comp}[\incl_2{\comp}\iota_r,\incl_2{\comp}C(\iota_{t})]^{\text{\rm rel}}
\\&\qquad\qquad\qquad
+ (((j^Q_1{\comp}\sigma(S^r))\times\sigma(S^r)){\comp}\Psi^2_0 \cup (\hat{\chi}_{\alpha}\vee\hat{\chi}_{\alpha})){\comp}[\incl_2{\comp}\iota_r,\incl_1{\comp}C(\iota_{t})]^{\text{\rm rel}}
\\&\qquad\qquad\qquad
+ (((j^Q_1{\comp}\sigma(S^r))\times\sigma(S^r)){\comp}\Psi^2_0 \cup (\hat{\chi}_{\alpha}\vee\hat{\chi}_{\alpha})){\comp} [\incl_1{\comp}\iota_r,\incl_2{\comp}C(\iota_{t})]^{\text{\rm rel}}
\\&\qquad
\sim \incl_1{\comp}(j^Q_1{\comp}\sigma(S^r){\comp}\Psi \cup \hat{\chi}_{\alpha}){\comp}[\iota_r,C(\iota_{t})]^{\text{\rm rel}}
\\&\qquad\qquad
+ \incl_2{\comp}(j^Q_1{\comp}\sigma(S^r){\comp}\Psi \cup \hat{\chi}_{\alpha}){\comp}[\iota_r,C(\iota_{t})]^{\text{\rm rel}}
\\&\qquad\qquad
+ [\hat{\chi}_{\alpha},j^Q_1{\comp}\sigma(S^r)]^{\text{\rm rel}}{\comp}\hat{T}
%\\&\qquad\qquad
+ [j^Q_1{\comp}\sigma(S^r),\hat{\chi}_{\alpha}]^{\text{\rm rel}},
\intertext{where $\hat{T} : S^{r+t}=S^{r-1}{\ast}S^{t} \to S^{t}{\ast}S^{r-1}=S^{r+t}$ is a switching map.
Since $[\hat{\chi}_{\alpha},j^Q_1{\comp}\sigma(S^r)]^{\text{\rm rel}} \sim \ast$ in $P^{\infty}(\G{Q}){\times}{\Sigma\G{Q}} \cup \{\ast\}{\times}P^{\infty}(\G{Q})$, we proceed as}&
\hat{\Delta}_{Q}{\comp}\sigma_0{\comp}\psi \sim \incl_1{\comp}\sigma_{\infty}{\comp}\psi + \incl_2{\comp}\sigma_{\infty}{\comp}\psi + [j^Q_1{\comp}\sigma(S^r),\hat{\chi}_{\alpha}]^{\text{\rm rel}}.
\intertext{On the other hand% in $P^{\infty}(\G{Q}) \vee P^{\infty}(\G{Q}) \cup {\Sigma}\G{Q}{\times}{\Sigma}\G{Q}$
, we have }&
\hat{\Delta}_{Q}{\comp}{\Sigma}\G{\psi}{\comp}\sigma(S^{r+t})
= (j^Q_1{\times}j^Q_1){\comp}\Delta_{{\Sigma}\G{Q}}{\comp}{\Sigma}\G{\psi}{\comp}\sigma(S^{r+t})
\\&\quad
= (j^Q_1{\times}j^Q_1){\comp}({\Sigma}\G{\psi}{\comp}\sigma(S^{r+t})\times{\Sigma}\G{\psi}{\comp}\sigma(S^{r+t})){\comp}\Delta_{S^{r+t}} 
\\&\quad
\sim (j^Q_1{\comp}{\Sigma}\G{\psi}{\comp}\sigma(S^{r+t}){\vee}j^Q_2{\comp}{\Sigma}\G{\psi}{\comp}\sigma(S^{r+t})){\comp}\mu_{r+t} 
\\&\quad
= \incl_1{\comp}\sigma_{\infty}{\comp}\psi + \incl_2{\comp}\sigma_{\infty}{\comp}\psi.
\end{align*}
Since a higher Hopf invariant $H^{\sigma_{0}}_2(\psi)$ is the difference between $\sigma_0{\comp}\psi$ and $j^Q_1{\comp}{\Sigma}\G{\psi}{\comp}\sigma(S^{r+t})$, we get the desired homotopy relation.
\end{Proof}
Next we show the following.
\begin{Prop}\label{prop:chi}~
There is a homotopy relation $\hat{\chi}_{\alpha} \sim j^Q_2{\comp}\chi_{p^{\G{Q}}_1}{\comp}C(H_1(\alpha))+j^Q_1{\comp}\delta_0 : (C(S^t),S^t) \to (P^{\infty}(\G{Q}),\im(j^Q_1{\comp}\sigma(S^r)))$ for some $\delta_0 : S^{t+1} \to {\Sigma}\G{Q}$, where the addition is given by the coaction $(C(S^{t}),S^t) \to (C(S^{t}){\vee}S^{t+1},S^t)$.
\end{Prop}
\begin{Proof}
Let $\chi'_{\alpha} : (C(S^{t}),S^t) \to (P^2(\G{Q}),\Sigma\G{Q})$ be the map given by the deformation of $\alpha$ to $p^{\G{Q}}_1{\comp}H_1(\alpha)$ in $\Sigma\G{Q}$ and by $\chi_{p^{\G{Q}}_1}{\comp}C(H_1(\alpha)) : (C(S^{t}),S^t) \to (P^2(\G{Q}),\Sigma\G{Q})$ as in \cite[Lemma 5.4, Remark 5.5]{Iwase:counter-ls}, where we denote by $C$ the functor taking cones.
Then by definition, we have $\chi'_{\alpha} \sim \chi_{p^{\G{Q}}_1}{\comp}C(H_1(\alpha))$ in $(P^2(\G{Q}),\Sigma\G{Q})$ and $j^Q_1{\comp}\chi'_{\alpha}\vert_{S^t}=j^Q_1{\comp}\sigma(S^r){\comp}\alpha=\hat{\chi}_{\alpha}\vert_{S^t}$.
%Then we obtain that $\sigma'_1{\comp}\chi_{\alpha} \sim \chi_{p^{\G{Q}}_1}{\comp}C(H_1(\alpha)) : (C(S^{t}),S^t) \to (P^2(\G{Q}),\Sigma\G{Q})$ by extending the above deformation.
%Let $\delta : S^{t+1} \to Q$ be 
Hence the difference between $\hat{\chi}_{\alpha}$ and $j^Q_2{\comp}\chi'_{\alpha}$ is given by a map $\delta : S^{t+1} \to P^{\infty}(\G{Q}){\simeq}Q$, which can be pulled back to $\delta_0 : S^{t+1} \to {\Sigma}\G{Q}$ ($\subset P^2(\G{Q})$) (see the proof of \cite[Theorem 5.6]{Iwase:counter-ls}).
%By $e^Q_2{\comp}\sigma_0{\sim}1_Q$, we have 
%the following factorisation
%\begin{equation*} 
%\chi_{\alpha} : (C(S^{t}),S^t) \overset{\chi_{\alpha}}\to
%(Q,S^r) \overset{\sigma_0}\to (P^2(\G{Q}),\Sigma{\G{S^r}}) \overset{e^Q_2}\to (Q,S^r),
%\end{equation*}
%and hence, 
Thus we have 
$%\begin{equation*}
\hat{\chi}_{\alpha} \sim
j^Q_2{\comp}\chi'_{\alpha}+j^Q_1{\comp}\delta_0 \sim
j^Q_2{\comp}\chi_{p^{\G{Q}}_1}{\comp}C(H_1(\alpha))+j^Q_1{\comp}\delta_0.
$%\end{equation*}
\end{Proof}
%as maps from $(C(S^{t}),S^t)$ to $(P^{\infty}(\G{Q}),\Sigma\G{Q})$, where the addition is given by the coaction $(C(S^{t}),S^t) \to (C(S^{t}){\vee}S^{t+1},S^t)$, and hence we have
Now we prove Lemma \ref{prop:standard_structure} using
Propositions \ref{prop:relative-whitehead} and \ref{prop:chi}.
\begin{align*}&
[j^Q_1,j^Q_2{\comp}\chi_{p^{\G{Q}}_1}]^{\text{\rm rel}}{\comp}H^{\sigma_{0}}_2(\psi) 
\sim \hat{\Delta}_Q{\comp}p^{\G{Q}}_2{\comp}H^{\sigma_{0}}_2(\psi) 
\sim [j^Q_1{\comp}\sigma(S^{r}),\hat{\chi}_{\alpha}]^{\text{\rm rel}} 
\\&\quad
\sim [j^Q_1{\comp}\sigma(S^r),j^Q_2{\comp}\chi_{p^{\G{Q}}_1}{\comp}C(H_1(\alpha))]^{\text{\rm rel}} + [j^Q_1{\comp}\sigma(S^r),j^Q_1{\comp}\delta_0]
\\&\quad
\sim \pm [j^Q_1,j^Q_2{\comp}\chi_{p^{\G{Q}}_1}]^{\text{\rm rel}}{\comp}(\hat{i}{\ast}1_{\G{Q}{\ast}\G{Q}}){\comp}(1_{S^{r-1}}{\ast}H_1(\alpha)) + (j^Q_1{\vee}j^Q_1){\comp}[\sigma(S^r),\delta_0] 
%\\
\intertext{by Propositions \ref{prop:relative-whitehead} and \ref{prop:chi}.
Since $[\sigma(S^r),\delta_0] \sim 0$ in $\Sigma\G{Q}{\times}\Sigma\G{Q}$, we proceed as}&
[j^Q_1,j^Q_2{\comp}\chi_{p^{\G{Q}}_1}]^{\text{\rm rel}}{\comp}H^{\sigma_{0}}_2(\psi) \sim \pm [j^Q_1,j^Q_2{\comp}\chi_{p^{\G{Q}}_1}]^{\text{\rm rel}}{\comp}(\hat{i}{\ast}1_{\G{Q}{\ast}\G{Q}}){\comp}{\Sigma}^{r}H_1(\alpha),
%\quad\text{since $[\sigma(S^r),\delta_0] \sim 0$ in $\Sigma\G{Q}{\times}\Sigma\G{Q}$}.
\end{align*}
%since $[\sigma(S^r),\delta_0] \sim 0$ in $\Sigma\G{Q}{\times}\Sigma\G{Q}$.
Since the relative Whitehead product $[j^Q_1,j^Q_2{\comp}\chi_{p^{\G{Q}}_1}]^{\text{\rm rel}}$ induces a split monomorphism in homotopy groups, we have $H^{\sigma_{0}}_2(\psi) \sim \pm (\hat{i}{\ast}1_{\G{Q}{\ast}\G{Q}}){\comp}{\Sigma}^{r}H_1(\alpha)$.
Thus we obtain 
$$
H^S_2(\psi) \ni [H^{\sigma_{0}}_2(\psi)] = \pm [(\hat{i}{\ast}1_{\G{Q}{\ast}\G{Q}}){\comp}{\Sigma}^{r}H_1(\alpha)].
$$
This completes the proof of Lemma \ref{prop:standard_structure}.
%\end{Proof}
%
%	MAIN RESULT
%
\section{Proof of Theorem \ref{lem:crucial}}\label{sect:crucial-lemma}
\par\noindent
Let $\beta : S^{v} \to S^{r+t}$ be a co-H-map.
If $[{\Sigma}^{r}H_1(\alpha){\comp}\beta] = 0$, then we have $H^S_2(\psi{\comp}\beta) \ni [H^{\sigma_{0}}_2(\psi){\comp}\beta] = \pm [(\hat{i}{\ast}1_{\G{Q}{\ast}\G{Q}}){\comp}{\Sigma}^{r}H_1(\alpha){\comp}\beta] = 0$ by Lemma \ref{prop:standard_structure}.
Hence we show the converse.
There are cofibre sequences as follows:
\begin{align*}&
S^{t} \overset{\alpha}\rightarrow S^r \overset{i}\hookrightarrow Q \overset{q}\rightarrow S^{t+1},\qquad
%\\&
S^{r+t} \overset{\psi}\rightarrow Q \overset{j}\hookrightarrow E \overset{\hat{q}}\rightarrow S^{r+t+1}.
\end{align*}
By the arguments in \S\ref{section:standard_structure}, we know there are `standard' structures $\sigma(S^r) : S^r \to P^1(\G{S^r})$ and $\sigma_{0} : Q \to P^2(\G{Q})$ for $\cat(S^r) = 1$ and $\cat(Q) = 2$, respectively, where $\sigma_{0}\vert_{S^{r}} = \sigma(S^{r})$ in $P^2(\G{Q})$.
%\begin{align*}&
%\sigma(S^r) : S^r \to P^1(\G{S^r})\quad \text{and}\quad 
%%\\&
%\sigma_{0} : Q \to P^2(\G{Q})\quad \text{($\sigma_{0}\vert_{S^{r}}=\sigma(S^{r})$ in $P^2(\G{Q})$)}
%\end{align*}
%where $\sigma_{0}$ is defined in Lemma \ref{prop:standard_structure}.
\par
Let $\sigma$ be a structure for $\cat(Q)=2$ with $H^{\sigma}_2(\psi){\comp}\beta \sim 0$ in $E^3(\G{Q})$.
For dimensional reasons, $\sigma{\vert}_{S^r}$ is homotopic to $\sigma(S^r)$ which is given by the bottom-cell inclusion.
We regard $e^{Q}_2 : P^2(\G{Q}) \to Q$ as a fibration with fibre $E^{3}(\G{Q}) \overset{p^{\G{Q}}_2}\to P^2(\G{Q})$ and $\sigma_{0}$ as a cross-section of $e^{Q}_2$.
Then by the definition of a structure, we have $e^{Q}_2{\comp}\sigma \sim 1_Q$.
Thus we obtain the following homotopy relations:
\begin{align*}&
\sigma\vert_{S^{r}} \sim \sigma(S^{r})=\sigma_{0}\vert_{S^{r}}\quad \text{in $P^2(\G{Q})$},\quad 
%\\&
e^{Q}_2{\comp}\sigma \sim e^{Q}_2{\comp}\sigma_{0}=1_Q.
\end{align*}
Thus the difference between $\sigma$ and $\sigma_{0}$ is given by a map $\gamma : S^{t+1} \to P^2(\G{Q})$ which can be lift to $E^3(\G{Q})$:
\begin{equation*}
%\sigma \sim \sigma' \sim \sigma'' \sim \sigma_{0} + p^{\G{Q}}_2{\comp}\gamma\quad \text{in $P^2(\G{Q})$,\quad $\gamma : S^{t+1} \to E^3(\G{Q})$},
\sigma \sim \sigma_{0} + \gamma\quad \text{in $P^2(\G{Q})$},
\end{equation*}
where the addition is taken by the coaction $\mu : Q \to Q \vee S^{t+1}$ along the collapsing $q : Q \to S^{t+1}$.
Thus we obtain that $\sigma{\comp}\psi \sim \{\sigma_{0},\gamma\}{\comp}\mu{\comp}\psi$ in $P^2(\G{Q})$, where $\{\sigma_{0},\gamma\} : Q \vee S^{t+1} \to P^{2}(\G{Q})$ is a map given by $\{\sigma_{0},\gamma\}\vert_{Q} = \sigma_{0}$ and $\{\sigma_{0},\gamma\}\vert_{S^{t+1}} = \gamma$.

By the definition of $\psi$, we have
\begin{equation*}
\mu{\comp}\psi \sim (\psi \vee q{\comp}\psi){\comp}\mu + [\iota'_{r},\iota''_{t+1}] \sim (\psi \vee \ast){\comp}\mu + [\iota'_{r},\iota''_{t+1}]\quad \text{in $Q \vee S^{t+1}$},
\end{equation*}
where $\iota'_r : S^r \hookrightarrow Q \hookrightarrow Q{\vee}S^{t+1}$ and $\iota''_{t+1} : S^{t+1} \hookrightarrow Q{\vee}S^{t+1}$ are inclusions.
Hence we have
\begin{equation*}
\sigma{\comp}\psi \sim \sigma_{0}{\comp}\psi + [\sigma(S^{r}),\gamma]\quad \text{in $P^2(\G{Q})$},
\end{equation*}
which yields the following homotopy relation, since $\beta$ is a co-H-map:
\begin{equation}\label{eq:hom_rel_Q1}
\begin{split}&
0 \sim p^{\G{Q}}_2{\comp}H^{\sigma}_2(\psi){\comp}\beta 
 \sim P^2(\G{\psi}){\comp}\sigma(S^{r+t}){\comp}\beta - \sigma{\comp}\psi{\comp}\beta
\\&
\sim P^2(\G{\psi}){\comp}\sigma(S^{r+t}){\comp}\beta - (\sigma_{0}{\comp}\psi{\comp}\beta + [\sigma(S^{r}),\gamma]{\comp}\beta) 
\\&
\sim (P^2(\G{\psi}){\comp}\sigma(S^{r+t}) - \sigma_{0}{\comp}\psi){\comp}\beta - [\sigma(S^{r}),\gamma]{\comp}\beta
\\&
\sim p^{\G{Q}}_2{\comp}H^{\sigma_{0}}_2(\psi){\comp}\beta - [\sigma(S^{r}),\gamma]{\comp}\beta 
\sim \pm p^{\G{S^{r}}}_2{\comp}{\Sigma}^rH_1(\alpha){\comp}\beta - [\sigma(S^{r}),\gamma]{\comp}\beta%\quad\text{in $P^2(\G{Q})$}.
\end{split}
\end{equation}
in $P^2(\G{Q})$.
To proceed further, we consider the following commutative ladder of fibre sequences.
\begin{equation*}
\begin{xy}
\xymatrix{
\G{S^{r}}~
	\ar@{^{(}->}[r] 
	\ar@{_{(}->}[d] 
&
E^3(\G{S^{r}}) 
	\ar[r]^{p^{\G{S^{r}}}_2} 
	\ar@{_{(}->}[d] 
&
P^2(\G{S^{r}}) 
	\ar[r]^{\quad e^{S^{r}}_2} 
	\ar@{_{(}->}[d] 
&
S^{r}
	\ar@{_{(}->}[d] 
\\
\G{Q}\, 
	\ar@{^{(}->}[r] 
&
E^3(\G{Q}) 
	\ar[r]^{p^{\G{Q}}_2} 
&
P^2(\G{Q}) 
	\ar[r]^{\quad e^{Q}_2} 
&
Q.
}
\end{xy}
\end{equation*}
Since the pair ($E^3(\G{Q})$,$E^3(\G{S^{r}})$) is ($t+2r-1$)-connected and $t+1 < r+t < t+2r-1$, $r > 1$, we have $\pi_{t+1}(E^3(\G{Q})) \cong \pi_{t+1}(E^3(\G{S^{r}}))$ and $\pi_{r+t}(E^3(\G{Q})) \cong \pi_{r+t}(E^3(\G{S^{r}}))$.
Since $\gamma$ can be lift to $E^3(\G{Q})$ and we know $\pi_{t+1}(E^3(\G{Q})) \cong \pi_{t+1}(E^3(\G{S^{r}}))$, we may regard that the image of $\gamma$ is contained in $P^2(\G{S^{r}})$.
Thus (\ref{eq:hom_rel_Q1}) implies a homotopy relation 
\begin{align}\label{eq:hom_rel_Q2}&
p^{\G{S^{r}}}_2\vert_{S^{r-1}{\ast}E^{2}(\G{S^r})}{\comp}{\Sigma}^rH_1(\alpha){\comp}\beta \sim \pm [\sigma(S^{r}),\gamma]{\comp}\beta\quad\text{in $P^2(\G{Q})$}.
\end{align}
Since $p^{\G{Q}}_2$ induces a split monomorphism in the homotopy groups of dimension $r+t$ and we know $\pi_{r+t}(E^3(\G{Q})) \cong \pi_{r+t}(E^3(\G{S^{r}}))$, (\ref{eq:hom_rel_Q2}) implies a homotopy relation 
\begin{equation*}%\label{eq:hom_rel_S}
p^{\G{S^{r}}}_2\vert_{S^{r-1}{\ast}E^{2}(\G{S^r})}{\comp}{\Sigma}^rH_1(\alpha) \sim \pm [\sigma(S^{r}),\gamma]\quad\text{in $P^2(\G{S^{r}})$}.
\end{equation*}

To show ${\Sigma}^rH_1(\alpha){\comp}\beta$ is trivial, we use the following proposition obtained by a straight-forward calculation (see Mac~Lane \cite{MacLane:homology}, Stasheff \cite{Stasheff:higher-associativity} or \cite{Iwase:k-ring}, for example) of Bar resolution:
\begin{Prop}~
The composition map 
$$
\partial : E^{m+1}(\G{S^{r}}) \overset{p^{\G{S^{r}}}_{m}}\to P^{m}(\G{S^{r}}) \to P^{m}(\G{S^{r}})/{\Sigma}\G{S^{r}} \simeq {\Sigma}E^{m}(\G{S^{r}})
$$ 
induces a homomorphism 
$$\partial_{\ast} : \tilde{H}_\ast(\wedge^{m+1}\G{S^r};\integral) \to \tilde{H}_\ast(\wedge^{m}\G{S^r};\integral),$$ 
which is given by 
\begin{align*}&
\partial_{\ast}(x^{a_0}{\otimes}x^{a_1}{\otimes}\cdots{\otimes}x^{a_m})
%\\&\qquad
= \sum^{m}_{i=1}(-1)^ix^{a_0}{\otimes}\cdots{\otimes}x^{a_{i-1}+a_{i}}{\otimes}\cdots{\otimes}x^{a_m},
\end{align*}
%\quad,
where $a_0,\cdots,a_m \geq 1$ and $x \in H_{r-1}(\G{S^{r}};\integral)$ is the generator of the Pontryagin ring $H_{\ast}(\G{S^{r}};\integral)$.
\end{Prop}
\begin{Cor}~
The composition map 
$$\partial' : S^{r-1}{\ast}E^{2}(\G{S^{r}}) \subset E^{3}(\G{S^{r}}) \overset{\partial}\to {\Sigma}E^{2}(\G{S^{r}}) \to {\Sigma}E^{2}(\G{S^{r}})/{\Sigma}(S^{r-1}{\ast}\G{S^{r}})
$$ 
induces an isomorphism 
$$\partial_{\ast} : \tilde{H}_\ast(S^{r-1}\wedge\G{S^r}\wedge\G{S^r};\integral) \to \tilde{H}_\ast((\G{S^r}/S^{r-1})\wedge\G{S^r};\integral),$$
which is given by 
\begin{equation*}
\partial'_{\ast}(x{\otimes}x^j{\otimes}x^k) = -x^{j+1}{\otimes}x^k,\quad j,k \geq 1.
\end{equation*}
\end{Cor}
Thus we obtain a left homotopy inverse of $p^{\G{S^{r}}}_2\vert_{S^{r-1}{\ast}E^{2}(\G{S^r})} : S^{r-1}{\ast}E^{2}(\G{S^{r}}) \to P^2(\G{S^{r}})$ as a composition map $P^2(\G{S^{r}}) \to P^2(\G{S^{r}})/{\Sigma}\G{S^{r}} \homeo {\Sigma}E^{2}(\G{S^{r}}) \to {\Sigma}E^{2}(\G{S^{r}})/{\Sigma}(S^{r-1}{\ast}\G{S^{r}}) \simeq S^{r-1}{\ast}E^{2}(\G{S^{r}})$, where the image of ${\Sigma}^rH_1(\alpha)$ lies in $S^{r-1}{\ast}E^{2}(\G{S^{r}})$.
On the other hand by the fact that $\im\sigma(S^{r}) \subset \Sigma\G{S^{r}}$, we also know that the Whitehead product $[\sigma(S^{r}),\gamma]$ vanishes in the quotient space $P^2(\G{S^{r}})/{\Sigma}\G{S^{r}}$, and hence never appears non-trivially in $S^{r-1}{\ast}E^2(\G{S^{r}})$.
Thus we conclude that ${\Sigma}^{r}H_1(\alpha){\comp}\beta$ is trivial.
%\begin{Rem}
%By using the homotopy-extension property for the pair ($Q$,$S^{r}$) and the homotopy-lifting-extension property for the fibration $e^{Q}_2$ with cross-section $\sigma_0$, we can easily see that $\gamma$ can be lift to a map $\gamma_0 : S^{t+1} \to E^3(\G{Q})$ to the fibre of $e^{Q}_2$.
%But it is not necessary for our case.
%%%
%By using the homotopy-extension property for the pair ($Q$,$S^{r}$), there is a map $\sigma'$ ($\sim \sigma$) with $\sigma'\vert_{S^{r}}=\sigma_{0}\vert_{S^{r}}$ and a homotopy $H : e^{Q}_2{\comp}\sigma' \sim e^{Q}_2{\comp}\sigma_{0}=1_Q$.
%Then $\sigma_{0}{\comp}H\vert_{S^{r}}$ gives a homotopy of $\sigma'\vert_{S^{r}}$ to $\sigma_{0}\vert_{S^{r}}$ in $P^2(\G{Q})$.
%By the homotopy-lifting-extension property for the fibration $e^{Q}_2$ and the pair ($Q$,$S^{r}$), there is a map $\sigma'' : Q \to P^2(\G{Q})$ and a lifting $\hat{H} : \sigma \sim \sigma''$ of the homotopy $H$ extending $\sigma_{0}{\comp}H\vert_{S^{r}}$.
%Hence we obtain that $\sigma''\vert_{S^{r}}=\sigma_{0}\vert_{S^{r}}$ and $e^{Q}_2{\comp}\sigma'' = e^{Q}_2{\comp}\sigma_{0}$.
%Thus the difference between $\sigma''$ and $\sigma_{0}$ is given by a map $\gamma$ to the fibre of $e^{Q}_2$:
%\end{Rem}
%
%	MAIN THEOREMS
%
%\section{Proofs of Theorems \ref{thm:main-theorem} and \ref{thm:cat-bundle_sphere}}\label{sect:main-theorem}
%\par\noindent
%By Theorem \ref{lem:crucial}, we obtain Theorem \ref{thm:main-theorem} or \ref{thm:cat-bundle_sphere} by just replacing $H^S_2(\psi)$ in Fact \ref{fact:cat-bundle1} or \ref{fact:cat-bundle2} with ${\Sigma}^{r}H_1(\alpha)$.

%
%     BIBLIOGRAPHY
%

%
\end{document}